\documentclass[12pt]{article}
\usepackage{amssymb,amsmath}
\usepackage{hyperref}
\usepackage{graphicx}
\usepackage{color}
\usepackage{chngcntr}
\usepackage[top=2.5cm,bottom=2cm,left=2.5cm,right=2.5cm]{geometry}

\begin{document}

\begin{center}
{\LARGE\bf Just another solution to the Basel Problem\\
\vskip 1cm}

\large
Alois Schiessl\\ 
\vskip 0.25cm
\tt{aloisschiessl@web.de}
\end{center}

\begin{abstract}
The Basel problem consists in finding the sum of the reciprocals of the squares of the positive integers. It was finally solved in 1735 by Leonhard Euler. He showed that 
\[
\sum_{k=1}^\infty \frac{1}{k^2}=\frac{\pi^2}{6}.
\]
In this paper, we propose a simple proof based on the Weierstrass Sine product formula and L'Hôpital's rule.
\end{abstract}
\centerline{\it In celebration of Pi Day 2023}
$ \\ $
$ \\ $
The Basel Problem was a very famous problem in the middle of the seventeenth century. It was first posed by Pietro Mengoli in 1650 and many prominent mathematicians of the time tried to solve it without success. It took almost a hundred years before Euler \cite{E1}, \cite{E2} succeeded in 1734 in proving the above closed-form solution.
$ \\ $
In this paper, we aim to give another proof of the Basel problem. The proof is short, simple and uses only classical analysis.
$ \\ \\ $
We recall the Weierstrass factorisation theorem \cite{RC} for $\sin \left(\pi\,x\right)$ : 
\begin{align}
\sin \left( {\pi x} \right) = \pi x \cdot \prod\limits_{k = 1}^\infty  {\left( {1 - \frac{{{x^2}}}{{{k^2}}}} \right)} 
\end{align}
\\
The infinite product is analytic, so we can take the natural logarithm on both sides:
\begin{align}
\ln \left( {\sin \left( {\pi x} \right)} \right) = \ln \left( {\pi x} \right) + \sum\limits_{k = 1}^\infty  {\ln \left( {1 - \frac{{{x^2}}}{{{k^2}}}} \right)}
\end{align}
Next differentiating with respect to $x$ and assuming $k^2-x^2\neq 0$ gives us
\begin{align}
\frac{{\pi\cdot  \cos \left( {\pi x} \right)}}{{\sin \left( {\pi x} \right)}}  = \frac{1}{x} + \sum\limits_{k = 1}^\infty  {\left( {\frac{1}{{1 - \frac{{{x^2}}}{{{k^2}}}}} \cdot \frac{{ - 2 x}}{{{k^2}}}} \right)}  = \frac{1}{x} - 2 x\sum\limits_{k = 1}^\infty  {\left( {\frac{1}{{{k^2} - {x^2}}}} \right)}
\end{align}
Uniform convergence allows the interchange of the derivative and the infinite series.
$ \\ $
Using elementary algebraic, equation (3)
can be rearranged to the equivalent form
\begin{align}
&\frac{{\pi\cdot  \cos \left( {\pi x} \right)}}{{\sin \left( {\pi x} \right)}}-\frac{1}{x}=- 2 x\sum\limits_{k = 1}^\infty  {\left( {\frac{1}{{{k^2} - {x^2}}}} \right)}\\
&\left(\frac{{\pi\cdot  \cos \left( {\pi x} \right)}}{{\sin \left( {\pi x} \right)}} -\frac{1}{x}\right)\cdot\frac{-1}{2\,x}
=\sum\limits_{k = 1}^\infty  {\left( {\frac{1}{{{k^2} - {x^2}}}} \right)}\\
&{\frac {\sin \left( \pi \,x \right) -\pi \,x\cos \left( \pi \,x
 \right) }{2\,{x}^{2}\sin \left( \pi \,x \right) }}
=\sum _{k=1}^{\infty }\frac {1}{ \left( k^2-x^2 \right)}
\,;\;\;k^2-x^2\neq 0
\end{align}
The study of equation (6) is central in the paper.
$ \\ $
$ \\ $
Note that if we take $x=0$, the right-hand side becomes the sum
\begin{equation}
\sum _{k=1}^
{\infty }\frac {1}{ k^2}=1+\frac {1}{2^2}+\frac {1}{3^2}+\frac {1}{4^2}+\ldots
\end{equation}
It is easy to verify that the series is convergent (Integral test for convergence).
$ \\ $
$ \\ $
We know that we must obtain something nice on the left-hand side because the right-hand side is nice. Plugging in $x=0$, the left-hand side unfortunately leads to an indeterminate form
\begin{align}
\lim_{x \to 0} \;\;{\frac {\sin \left( \pi \,x \right) -\pi \,x\cos \left( \pi \,x
 \right) }{2\,{x}^{2}\sin \left( \pi \,x \right) }}
\rightarrow \frac{0}{0}
\end{align}
\\From the convergence of the right hand side, we conclude that
\begin{align}
\lim_{x \to 0} \;\;{\frac {\sin \left( \pi \,x \right) -\pi \,x\cos \left( \pi \,x
 \right) }{2\,{x}^{2}\sin \left( \pi \,x \right) }}
\end{align}
exists. We use L'Hôpital's rule to evaluate the limit. Applying L'Hôpital's rule once still gives an indeterminate form. In this case, the limit can be evaluated by applying the rule twice more.
$ \\ $
$ \\ $
$1.$ Application of L'Hôpital's rule:
\begin{align}
&\lim_{x \to 0} \;\;{\frac {\sin \left( \pi \,x \right) -\pi \,x\cos \left( \pi \,x
 \right) }{2\,{x}^{2}\sin \left( \pi \,x \right) }}\\
=&\lim_{x \to 0} \;{\frac {{\pi }^{2} x\sin \left( \pi \,x \right)}{4\,x\sin \left( \pi \,x \right) +2\,\pi \,{x}^{2}\cos \left( \pi \,x
 \right)}}\rightarrow\;\frac{0}{0}
\end{align}
$2.$ Application of L'Hôpital's rule:
\begin{align}
&\lim_{x \to 0} \;\;{\frac {\sin \left( \pi \,x \right) -\pi \,x\cos \left( \pi \,x
 \right) }{2\,{x}^{2}\sin \left( \pi \,x \right) }}\\
=&\lim_{x \to 0} \frac {{\pi }^{3}\cos \left( \pi \,x \right) x+{\pi }^{2}\sin \left( \pi \,x
 \right)}{4\,\sin \left( \pi \,x \right) +8\,\pi \,\cos \left( \pi \,x \right) x
-2\,{\pi }^{2}{x}^{2}\sin \left( \pi \,x \right)}\rightarrow\;\frac{0}{0}
\end{align}
$3.$ Application of L'Hôpital's rule:
\begin{align}
&\lim_{x \to 0} \;\;{\frac {\sin \left( \pi \,x \right) -\pi \,x\cos \left( \pi \,x
 \right) }{2\,{x}^{2}\sin \left( \pi \,x \right) }}\\
=\;&\lim_{x \to 0} \;{\frac {2\,{\pi }^{3}\cos \left( \pi 
\,x \right)-{\pi }^{4}\sin \left( \pi \,x \right) x 
}{12\,\pi \,\cos \left( \pi \,x \right) -12\,{\pi }^{2}\sin \left( \pi 
\,x \right) x-2\,{\pi }^{3}{x}^{2}\cos \left( \pi \,x \right) }}\\
=\;&\frac {2\,{\pi }^{3}\cos \left( \pi 
\,\cdot 0\right)-{\pi }^{4}\sin \left( \pi \,\cdot 0 \right) \cdot 0 
}{12\,\pi \,\cos \left( \pi \,\cdot 0 \right) -12\,{\pi }^{2}\sin \left( \pi 
\,\cdot 0 \right) \cdot {0}^{2}\cdot\,{\pi }^{3}\cdot {0}^{2}\cdot\cos \left( \pi \,\cdot 0 \right) }\\
=\;&\frac {2\,\pi^3}{12\,\pi}=\frac {\pi^2}{6}
\end{align}
\\Putting all the pieces together, we immediately obtain the famous Euler identity
\[
\sum _{k=1}^
{\infty }\frac { 1}{ k^2}\,=\,\frac {\pi^2}{6}.
\]

\end{document}